\documentclass[12pt]{article}
\usepackage{amssymb,amsmath,amsfonts,latexsym,amsthm,geometry}
\usepackage{amsmath}
\usepackage{amsfonts}
\usepackage{amssymb}
\usepackage{graphicx}
\usepackage{hyperref}
\providecommand{\U}[1]{\protect\rule{.1in}{.1in}}
\providecommand{\U}[1]{\protect\rule{.1in}{.1in}}
\providecommand{\U}[1]{\protect\rule{.1in}{.1in}}
\newtheorem{theorem}{Theorem}[section]
\newtheorem{corollary}[theorem]{Corollary}
\newtheorem{proposition}[theorem]{Proposition}
\newtheorem{lemma}[theorem]{Lemma}

\newtheorem{Fact}[theorem]{Fact}

\theoremstyle{definition}

\newtheorem{remark}[theorem]{Remark}

\begin{document}

\title{Some Notes on Complex Symmetric Operators}

\author{\sc Marcos S. Ferreira
\thanks{DCET, UESC, Bahia, Brazil, msferreira@uesc.br}}

\maketitle

\begin{abstract}
In this paper we show that every conjugation $C$ on the Hardy-Hilbert space $H^{2}$ is of type $C=T^{*}C_{1}T$, where $T$ is an unitary operator and $C_{1}f\left(z\right)=\overline{f\left(\overline{z}\right)}$, with $f\in H^{2}$. In the sequence, we extend this result for all separable Hilbert space $\mathcal H$ and we prove some properties of complex symmetry on $\mathcal H$. Finally, we prove some relations of complex symmetry between the operators $T$ and $\left|T\right|$, where $T =U\left|T\right|$ is the polar decomposition of bounded operator $T\in\mathcal L\left(\mathcal H\right)$ on the separable Hilbert space $\mathcal H$.
\end{abstract}

\section{Introduction}

Let $\mathcal L\left(\mathcal H\right) $ be the space of bounded linear operators on a
separable Hilbert space $\mathcal H$. \ A \textit{conjugation} $C$ on $\mathcal H$ is an
antilinear operator $C:\mathcal H\rightarrow\mathcal H$ such that $C^{2}=I$ and $\left\langle
Cf,Cg\right\rangle =\left\langle g,f\right\rangle ,$ for all $f,g\in\mathcal H$. An operator $T\in \mathcal L\left( H\right) $ is
said to be \textit{complex symmetric} if there exists a conjugation $C$ on $H
$ such that $CT=T^{\ast }C$ (we will often say that $T$ is $C$-symmetric). Complex symmetric operators generalize the concept of symmetric matrices of linear algebra. Indeed, it is well known (\cite[Lemma 1]{F}) that given a conjugation $C$, there exists an orthonormal basis $\left\{
f_{n}\right\} _{n=0}^{\infty }$ for $\mathcal H$ such that $Cf_{n}=f_{n}$. Hence, if $T$ is $C$-symmetric then
\begin{eqnarray}\label{eq1}
\left\langle Tf_{n},f_{m}\right\rangle=\left\langle Cf_{m},CTf_{n}\right\rangle=\left\langle f_{m},T^{*}Cf_{n}\right\rangle=\left\langle Tf_{m},f_{n}\right\rangle,
\end{eqnarray}
that is, $T$ has a symmetric matrix representation. The reciprocal of this fact is also true. That is, if there is an orthonormal basis such that $T$ has a symmetric matrix representation, then $T$ is complex symmetric. 

The complex symmetric operators class was initially addressed by Garcia and Putinar \cite{F,B} and includes the normal operators, Hankel operators and Volterra integration operators.

Now, let $L^{2}$ be the Hilbert space on the unit circle $\mathbb{T}$ and let $%
L^{\infty }$ be the Banach space of all essentially bounded functions on $\mathbb{T}$. It is known that $\left\{ e_{n}\left( e^{it}\right) :=e^{int}:n\in\mathbb{Z}\right\} $ is an orthonormal basis for $L^{2}$. The \textit{Hardy-Hilbert space}, denoted by $H^{2}$, consists of all analytic functions $f\left( z\right)
=\sum_{n=0}^{\infty }a_{n}z^{n}$ on the unit disk $\mathbb D$ such that $\sum_{n=0}^{\infty }\left\vert
a_{n}\right\vert ^{2}<\infty .$ It is clear that $\mathcal B:=\left\{e_{n}(z)=
z^{n}:n=0,1,2,\ldots \right\} $ is an orthonormal basis for $H^{2}.$

For each $\phi \in L^{\infty },$ the \textit{Toeplitz operator} $T_{\phi
}:H^{2}\rightarrow H^{2}$ is defined by%
\begin{equation*}
T_{\phi }f=P\left( \phi f\right) ,
\end{equation*}%
for each $f\in H^{2},$ where $P:L^{2}\rightarrow H^{2}$ is the orthogonal
projection. The concept of Toeplitz operators was initiated by Brown and Halmos \cite{A} and generalizes the concept of Toeplitz matrices.

In \cite{C}, Guo and Zhu raised the question of characterizing complex symmetric Toeplitz operators on $H^{2}$ in the unit disk. In order to obtain such characterization, Ko
and Lee \cite{D} introduced the family of conjugations $C_{\lambda}:H^{2}\rightarrow H^{2},$ given by 
\begin{equation*}
C_{\lambda }f\left( z\right) =\overline{f\left( \lambda \overline{z}%
\right) }
\end{equation*}%
with $\lambda \in\mathbb{T}$ and proved the following result:

\begin{theorem}\label{teo3}
Let $\phi(z)=\sum_{n=-\infty}^{\infty}\widehat{\phi}(n)z^{n}\in L^{\infty}$. Then $T_{\phi}$ is $C_{\lambda}$-symmetric if, and only if, $\widehat{\phi}(-n)=\lambda^{n}\widehat{\phi}(n)$, for all $n\in\mathbb Z$. 
\end{theorem}

\section{Canonical Conjugations}

Our first objective in this paper is to study relations between an arbitrary conjugation $C$ on $H^{2}$ and the conjugation $C_{1}f\left(z\right)=\overline{f\left(\overline{z}\right)}$. Once the conjugation $C_{1}$ is a kind of canonical conjugation on $H^{2}$, we observe a close relationship between conjugations of $H^{2}$ and conjugation $C_{1}$, namely:

\begin{theorem}\label{teo1}
If $C$ is an conjugation on $H^{2},$ then exists an unitary operator $%
T:H^{2}\rightarrow H^{2}$ such that $TC=C_{1}T.$
\end{theorem}
\begin{proof}
Since $C$ is an conjugation, there exists an orthonormal
basis $\mathcal B^{\prime}=\left\{ f_{n}\right\} _{n=0}^{\infty }$ of $H^{2}$
such that $Cf_{n}=f_{n}.$ Now, let $\mathcal B=\left\{e_{n}\right\}_{n=0}^{\infty}$ the standard orthonormal basis of $H^{2}$ and the linear isomorphism $T:H^{2}\rightarrow H^{2}$ given by
\begin{equation*}
T\left(\sum\limits_{n=0}^{\infty
}a_{n}f_{n}\right) =\sum\limits_{n=0}^{\infty }a_{n}e_{n}.
\end{equation*}

Note that $Tf_{n}=e_{n}$, for all $n\geq 0$, and therefore $T$ is unitary. Now, for $f\left(z\right)=\sum_{n=0}^{\infty }a_{n}e_{n}\in H^{2}$, we get
\begin{eqnarray*}
C_{1}f\left( z\right) &=&\sum\limits_{n=0}^{\infty }\overline{a_{n}}e_{n} \\
&=&\sum\limits_{n=0}^{\infty }\overline{a_{n}}T\left( f_{n}\right) \\
&=&T\left( \sum\limits_{n=0}^{\infty }\overline{a_{n}}Cf_{n}\right) \\
&=&\left( TC\right) \left( \sum\limits_{n=0}^{\infty
}a_{n}f_{n}\right) \\
&=&\left( TC\right) \left( \sum\limits_{n=0}^{\infty
}a_{n}T^{-1}\left( e_{n}\right) \right) \\
&=&\left( TCT^{-1}\right) f\left( z\right) ,
\end{eqnarray*}%
whence $C_{1}T=TC.$
\end{proof}

The previous theorem says that all complex symmetric Toeplitz operator is unitarily equivalent to a $C_{1}$-symmetric operator. Indeed:

\begin{remark}
Let $T_{\phi }:H^{2}\rightarrow H^{2}$ an Toeplitz operator. Observe that, if $T_{\phi }$ is $C$-symmetric, since the
operator $T$ of previous theorem is unitary, we have%
\begin{equation*}
C_{1}=TCT^{\ast },
\end{equation*}%
therefore the operator $%
T_{2}:=TT_{\phi }T^{\ast }$ is $C_{1}$-symmetric (see \cite[p. 1291]{F}). This shows that $T_{\phi }$
and $T_{2}$ are unitarily equivalent operators. Moreover, is obvious that,
if $T$ commutes with $C_{1}$ or $C$, then $C=C_{1}.$
\end{remark}

\begin{corollary}
Let $A\in \mathcal L\left( H^{2}\right) .$ Then $A$ is $C_{1}$-symmetric if, and only
if, the matrix of $A$ with respect the canonical basis of $H^{2}$ is
symmetric.
\end{corollary}

\begin{proof}
If $A$ is $C_{1}$-symmetric, then $C_{1}A=A^{\ast }C_{1}.$
Moreover, by previous theorem there exists an isomorphism $T$ on $H^{2}$
such that $TC_{1}=C_{1}T.$ Consider $\mathcal B=\left\{ e_{n}\right\}
_{n=0}^{\infty }$ and $\mathcal B^{\prime }=\left\{ f_{n}\right\} _{n=0}^{\infty
}$ orthonormal basis of $H^{2}$ such that%
\begin{equation*}
Tf_{n}=e_{n}\text{ and }C_{1}f_{n}=f_{n}.
\end{equation*}

Thus, we must%
\begin{equation*}
C_{1}e_{n}=C_{1}\left( Tf_{n}\right) =TC_{1}\left( f_{n}\right)
=Tf_{n}=e_{n},
\end{equation*}%
that is $C_{1}e_{n}=e_{n},$ $\forall n\geq 0.$ Therefore, by \eqref{eq1}, follows that $\left[ A\right] _{\mathcal B}=\left[ A\right] _{\mathcal B}^{t}.$

Reciprocally, suppose that $A$ is $C$-symmetric such that $Ce_{n}=e_{n}.$ By
previous theorem, $TC=C_{1}T$ and $Te_{n}=e_{n}.$ Hence, $T$ is the identity
operator and so $C=C_{1}.$
\end{proof}

In fact, the reciprocal of the Theorem \ref{teo1} is true:

\begin{proposition}\label{prop1}
If $T:H^{2}\rightarrow H^{2}$ is an unitary operator, then $C:=T^{-1}C_{1}T$ is an conjugation on $H^{2}$.
\end{proposition}

\begin{proof}
It is easy to see that $C$ is an antilinear operator. Now, since $T$ is an unitary operator, considering $\mathcal B=\left\{ e_{n}\right\}
_{n=0}^{\infty }$ the orthonormal basis of $H^{2}$, we have
\begin{eqnarray*}
\left\langle Ce_{n},Ce_{m}\right\rangle &=&\left\langle T^{*}C_{1}Te_{n},T^{*}C_{1}Te_{m}\right\rangle \\
&=&\left\langle TT^{*}C_{1}Te_{n},C_{1}Te_{m}\right\rangle \\
&=&\left\langle Te_{m},Te_{n}\right\rangle \\
&=&\left\langle e_{m},T^{*}Te_{n}\right\rangle \\
&=& \left\langle e_{m},e_{n}\right\rangle.
\end{eqnarray*}
By other hand, once $C^{2}=\left(T^{-1}C_{1}T\right)\left(T^{-1}C_{1}T\right)=I$, follow the desired.
\end{proof}

In short, the Theorem \ref{teo1} and the Proposition \ref{prop1} tell us that:

\begin{corollary}\label{cor1}
If $T:H^{2}\rightarrow H^{2}$ an linear isomorphism and $C:=T^{-1}C_{1}T$, then $T$ is unitary if, and only if, $C$ is a conjugation on $H^{2}$.
\end{corollary}

Now, once every separable Hilbert space has an orthonormal basis, follows that the Corollary \ref{cor1} is true for any separable Hilbert space $\mathcal H$. In fact, if $\mathcal B=\left\{f_{n}\right\}$ is an orthonormal basis on $\mathcal H$, then $\mathcal J:\mathcal H\rightarrow\mathcal H$ given by
\begin{eqnarray}\label{eq2}
\mathcal J\left(\sum\limits_{n=0}^{\infty
}\lambda_{n}f_{n}\right) =\sum\limits_{n=0}^{\infty }\overline{\lambda_{n}}f_{n}.
\end{eqnarray}
is a conjugation on $\mathcal H$. Thus, we have:

\begin{theorem}
If $T:\mathcal H\rightarrow\mathcal H$ an linear isomorphism and $C:=T^{-1}\mathcal J T$, then $T$ is unitary if, and only if, $C$ is a conjugation on $\mathcal H$.
\end{theorem}
\begin{proof}
Analogous to Theorem \ref{teo1} and Proposition \ref{prop1}.
\end{proof}

\begin{remark}
Note that in the Hardy-Hilbert space $H^{2}$, we have $\mathcal J=C_{1}$. 
\end{remark}

We already know that every normal operator is complex symmetric and that the reciprocal in general is not true. However, for Toeplitz operators, Theorem \ref{teo3} gives us:

\begin{Fact}
If $T_{\phi}$ is $\mathcal J$-symmetric, then $T_{\phi}$ is normal.
\end{Fact}

Now note that if $T_{\phi}$ is normal not necessarily $T_{\phi}$ is $\mathcal J$-symmetric. In fact, if $\phi(z)=-\overline{z}+z$ then $T_{\phi}$ is normal, however is not $\mathcal J$-symmetric.

\section{Properties of Complex Symmetry}

In the following, we present some properties of complex symmetry in Hilbert spaces. The first result gives us a way to get complex symmetric operators from another complex symmetric operator. First, we need some lemmas:

\begin{lemma}(\cite[Lemma 1]{B})\label{lem1}
If $C$ and $J$ are conjugations on a Hilbert space $\mathcal H$, then $U=CJ$ is a unitary operator. Moreover, $U$ is both $C$-symmetric and $J$-symmetric.
\end{lemma}

\begin{lemma}(\cite[Lemma 2.2]{G})\label{lem2}
If $U:\mathcal H\rightarrow\mathcal H$ is a unitary and complex symmetric operator with conjugation $C$, then $UC$ is a conjugation.
\end{lemma}

\begin{proposition}\label{prop2}
Let $T:\mathcal H\rightarrow\mathcal H$ an operator and $C$ and $J$ conjugations on $\mathcal H$. Then $T$ is $C$-symmetric if, and only if, $UT$ is $UC$-symmetric, where  $U=CJ$. 
\end{proposition}
\begin{proof}
We already know that $U$ is unitary and $C$ and $J$-symmetric and that $UC=CJC$ is a conjugation, by Lemmas \ref{lem1} and \ref{lem2}. Now since $U^{*}=U^{-1}=JC$ and $T$ is $C$-symmetric, we have
$$
UT\left(UC\right)=UTCU^{*}=UCT^{*}U^{*}=UC(UT)^{*}.
$$
Reciprocally, suppose that $UC(UT)^{*}=UT(UC)$. Thus
\begin{eqnarray*}
CT^{*}U^{*}&=&C(UT)^{*} \\
&=&U^{*}UC(UT)^{*} \\
&=&U^{*}UTUC \\
&=&TUC \\
&=&TCU^{*},
\end{eqnarray*}%
whence $CT^{*}=TC$.
\end{proof}

\begin{lemma}\label{lem3}
If $T:\mathcal H\rightarrow\mathcal H$ is both $C$-symmetric and $J$-symmetric, then $T$ is both $CJC$-symmetric and $JCJ$-symmetric.
\end{lemma}
\begin{proof}
By Lemma \ref{lem1}, we have that $U:=CJ$ is unitary and $C$ and $J$-symmetric. Hence, by Lemma \ref{lem2}, $UC=CJC$ is a conjugation on $\mathcal H$. Thus, since $CT=T^{*}C$ and $JT=T^{*}J$ we get
$$
\left(CJC\right)T=C\left(TJ\right)C=T^{*}\left(CJC\right),
$$
and so $T$ is $CJC$-symmetric. Analogous, we prove that $T$ is $JCJ$-symmetric. 
\end{proof}

\begin{proposition}
If $T:\mathcal H\rightarrow\mathcal H$ is both $C$ and $J$-symmetric, then $TU$ is $C$-symmetric, where $U=CJ$.
\end{proposition}
\begin{proof}
In fact, once $T$ is both $C$-symmetric and $J$-symmetric, we have by Lemma \ref{lem3} that $T$ is $CJC$-symmetric and so
$$
\left(TU\right)C=T\left(CJC\right)=CU^{*}T^{*}=C\left(TU\right)^{*}.
$$
\end{proof}

\begin{proposition}
Let $A:\mathcal H\rightarrow\mathcal H$ an invertible operator and $C$-symmetric. If $T$ is an operator on $\mathcal H$ such that $TA=AT$, then $T$ is $C$-symmetric if, and only if, $TA$ is $C$-symmetric.
\end{proposition}

\begin{proposition}
Let $U:\mathcal H\rightarrow\mathcal H$ an unitary operator $J$-symmetric. If $T$ is an operator such that $UT^{*}=TU$ (that is, $T$ and $T^{*}$ are unitarily equivalents), then:
\\ (i) $JT^{*}=T^{*}J\Leftrightarrow T$ is \ $UJ-$symmetric.
\\ (ii) $UJT=TJU^{*}\Leftrightarrow T$ is $J-$symmetric.
\end{proposition}

\begin{proposition}
An operator $T:\mathcal H\rightarrow\mathcal H$ is $C$-symmetric if, and only if, $\mathcal J T^{*}C=(C\mathcal J)^{*}T$.
\end{proposition}
\begin{proof}
We already know that $U=C\mathcal J$ is unitary and both $C$ and $\mathcal J$-symmetric. Now, note that
$$
\mathcal J T^{*}C=(C\mathcal J)^{*}T\Leftrightarrow UT^{*}C=CU^{*}T.
$$ 
First see that if $T$ is $C$-symmetric, then $UT^{*}C=U(CT)=(CU^{*})T$. Reciprocally, we have
\begin{eqnarray*}
CT^{*} &=&CU^{*}(UT^{*}C)C \\
&=&CU^{*}(CU^{*}T)C\\
&=&(UCCU^{*})TC \\
&=&TC.
\end{eqnarray*}%
\end{proof}

\begin{proposition}
Let $T:\mathcal H\rightarrow\mathcal H$ an operator and $C$ a conjugation on $\mathcal H$. If $TC=CT$, then $T$ is $C$-symmetric if, and only if, $T$ is self-adjoint.
\end{proposition}

\section{Complex Symmetry of Aluthge and Duggal Transforms}

Recall that the polar decomposition of an operator $T:\mathcal H\rightarrow\mathcal H$ is uniquely expressed by $T=U\left|T\right|$, where $\left|T\right|=\sqrt{T^{*}T}$ is a positive operator and $U$ is a partial isometry such that $Ker(U)=Ker\left|U\right|$ and $U$ maps cl($Ran\left|T\right|$) onto cl($Ran(T)$). In this case, the Aluthge and Duggal Transforms are given, respectively, by $\widetilde{T}=\left|T\right|^{\frac{1}{2}}U\left|T\right|^{\frac{1}{2}}$ and $\widehat{T}=\left|T\right|U$.

We already known that the Aluthge transform of a complex symmetric operator is also complex symmetric (see \cite[Theorem 1]{H}). In this section we study relations between complex symmetry of $T$ and $\left|T\right|$ with relation the conjugations $C$ and $J$, as well as the operators $\widetilde{T}$ and $\widehat{T}$.

\begin{proposition}
If $T$ is complex symmetric, then $\left|T\right|$ is also complex symmetric.
\end{proposition}
\begin{proof}
If $CT=T^{*}C$, we have by Remark of \cite[Lemma 1]{H} that $T=CJ\left|T\right|$, where $J$ commutes with $\left|T\right|$. Thus, once that $CJ$ is a unitary operator, follows that
$$
J\left|T\right|=C(CJ\left|T\right|)=\left|T\right|^{*}(CJ)^{*}C=\left|T\right|^{*}J.
$$
\end{proof}

\begin{corollary}
If $T$ is complex symmetric, then $\left|T\right|$ is self-adjoint.
\end{corollary}

\begin{proposition}
Let $C$ and $J$ conjugations on $\mathcal H$ such that $T=CJ\left|T\right|$. If $\left|T\right|$ is $C$-symmetric, then $T$ is also $C$-symmetric. 
\end{proposition}
\begin{proof}
First, let's show that $\left|T\right|$ is $J$-symmetric. In fact, see that
$$
J(JC\left|T\right|)=C\left|T\right|=\left|T\right|^{*}C=(\left|T\right|^{*}CJ)J,
$$
and so $JC\left|T\right|$ is $J$-symmetric. Thus, by Proposition \ref{prop2}, $\left|T\right|$ is $J$-symmetric.
Therefore, it is enough to see that:
\begin{eqnarray*}
CT &=&C(CJ\left|T\right|) \\
&=&\left|T\right|^{*}J \\
&=&(\left|T\right|^{*}JC)C \\
&=&(CJ\left|T\right|)^{*}C \\
&=& T^{*}C.
\end{eqnarray*}%
\end{proof}

\begin{corollary}
Let $T=CJ\left|T\right|$. If $\left|T\right|$ is $C$-symmetric, then $\widehat{T}=T$.
\end{corollary}

\begin{corollary}
Let $T=CJ\left|T\right|$. Then $\left|T\right|$ is $C$-symmetric if, and only if, $\widehat{T}$ is $J$-symmetric.
\end{corollary}

\begin{proposition}
Let $T=CJ\left|T\right|$. If $C\left|T\right|=\left|T\right|^{*}C$ and $CJ=JC$, then $T$ is $J$-symmetric.
\end{proposition}
\begin{proof}
In fact, we have that
\begin{eqnarray*}
JT &=&J(CJ\left|T\right|) \\
&=&C\left|T\right| \\
&=&\left|T\right|^{*}JJC \\
&=&\left|T\right|^{*}JCJ \\
&=& (CJ\left|T\right|)^{*}J \\
&=& T^{*}J.
\end{eqnarray*}%
\end{proof}

\end{document}